\documentclass{elsart}
\usepackage{latexsym}
\usepackage{amssymb}
\usepackage{amsfonts}
\usepackage{amsmath}
\usepackage{graphicx}
\usepackage{natbib}
\journal{Topology and its Applications}

\begin{document}

\begin{frontmatter}

\title{Polygonal approximation and energy of smooth knots}
\author{Eric J.~Rawdon\corauthref{cor}\thanksref{eric}}
\corauth[cor]{Corresponding author.}
\ead{rawdon@mathcs.duq.edu}
\address{Department of Mathematics and Computer Science\\Duquesne
University\\Pittsburgh, PA 15282}
\ead[url]{http://www.mathcs.duq.edu/\textasciitilde{rawdon}}
\thanks[eric]{This material is based upon work supported by the National Science Foundation under Grant No.~0074315 and 0296098.}

\author{Jonathan K.~Simon\thanksref{jon}}
\address{Department of Mathematics\\ University of Iowa\\
Iowa City, IA 52242}
\ead{jsimon@math.uiowa.edu}
\ead[url]{http://www.math.uiowa.edu/\textasciitilde{jsimon}}
\thanks[jon]{This material is based upon work supported by the National Science Foundation under Grant No.~9706789 and 0107209.}

\begin{abstract}
We establish a fundamental connection between smooth and polygonal
knot energies, showing that the \textit{Minimum Distance Energy} for
polygons inscribed in a smooth knot converges to the \textit{M\"obius
Energy} of the smooth knot as the polygons converge to the smooth
knot.  However, the polygons must converge in a ``nice'' way, and the 
energies must be correctly regularized.
We determine an explicit error bound between the energies
in terms of the number of the edges of the polygon and the
\textit{Ropelength} of the smooth curve.
\end{abstract}


\end{frontmatter}

\def\K{$K$}
\def\d{\,{\mathrm{d}}}
\def\arc#1{{\hbox{$\mathrm{arc}(#1)$}}}
\def\chord#1{\hbox{{$\mathrm{chord}(#1)$}}}
\def\umd{U_{md}}
\def\upmd{U'_{md}}
\def\emd{E_{md}}
\def\sbo{\hbox{$\mathrm{summand}$}}
\def\teb{\hbox{$\mathrm{total\ error}$}}
\def\md{md}
\def\elk{E_L(K)}

\section{Introduction} 

Given a knot $K$ in 3-space, there are
several ways to define an ``energy function'' that measures how
complicated the knot is in its spatial conformation.  In this paper,
we establish a fundamental approximation theorem, showing that when
both are appropriately normalized, the \textit{Minimum Distance
Energy} for polygons inscribed in a smooth curve converge to the
\textit{M\"obius Energy} of the curve as the polygons converge
to the smooth knot.  We do a careful analysis, and determine an
explicit error bound (Theorem \ref{explicit}), from which the
approximation (Theorem \ref{approaches}) follows immediately.

In  Section \ref{MainResultsAndNotation}, we state the main
theorems and agree on notation for the whole paper.  In Section
\ref{lemmas} we present a number of lemmas.  These establish
useful properties of curves and chords, so they may be of
independent interest. In Section
\ref{outlineproof}, we outline the proof of the error bound
(Theorem
\ref{explicit}), especially how to divide the problem into
several cases (more precisely, divide the domains into
different ``zones'') for which different analyses are
needed.   In Section
\ref{proofsforzones},  we give the detailed analyses for the
various cases, and in Section
\ref{finalproof}, we combine the results from Section
\ref{proofsforzones} to obtain the overall bound.

Of course the error depends on how well the polygon approximates
the smooth curve.
However, there are more subtle
issues to confront in controlling the error: One must reckon
with the amount of curvature the knot has, and how close it is
to being self-intersecting.  These are captured by the {\it
thickness radius}
$r(K)$ (see later in this section for definition).
Our error bound is developed
in terms of  the total arc-length
$\ell(K)$, the number of edges of the inscribed polygon $n$,
the mesh size
$\delta =\frac{\ell(K)}{n}$, the thickness radius $r(K)$, and the ratio
$E_L(K) =\frac{\ell(K)}{r(K)}$.  
Since these quantities are interrelated, there are
various ways to write the bound: the one we give in Theorem
\ref{explicit} is stated in terms of $n$ and $E_L(K)$ to
emphasize that it is invariant under change of scale.

Let $t \rightarrow x(t)$ be a unit-speed parameterization of
$K$ with domain a circle C. The {\it M\"obius Energy} or {\it
O'Hara Energy} is
$$E_0(K) = \iint_{C \times C}
\;\;\;\frac{1}{|x(t)-y(s)|^2}-
\frac{1}{|s-t|^2}\;\;ds\;dt. $$

The energy $E_{0}$ was defined and studied in
\cite{O1,FHW,KS,kusnerkim}.
The subscript $0$ in
$E_0$ reminds us that this version of the energy is exactly
zero if $K$ is a circle.

By visualizing a smooth knot as being 
made  of some  ``rope'', with a positive thickness,
we obtain a fundamental measure of knot complexity.
Hold the
core knot $K$ fixed and thicken the rope until the moment of
self-contact. Call that sup radius the {\it thickness radius}
or {\it injectivity radius},
$r(K)$.  Here is a more precise definition:  For small enough
$r$, the knot $K$ has a solid torus neighborhood consisting of pairwise
disjoint disks of radius
$r$ centered at the points of $K$ and orthogonal to $K$ at
those centers. Gradually increase $r$ until
some meridional disks touch; we call that supremum of
good radii $r(K)$.
The ratio 
$$E_L(K) = \frac{\text{arc-length of } K}{r(K)}\;,$$
called the {\it Rope-Length} of $K$, is a scale-invariant
numerical measure of knot compaction.

The basic theorems on thickness appear in
\cite{LSDR}, although the energy $E_L$ first appeared in
\cite{BO1}.  We recall the properties of $E_L$ in Section
\ref{lemmasaboutthickness}.

Let $P$ be a polygon with $n$ edges. The {\it Minimum Distance
Energy} 
 of $P$ is defined \cite{Si2} as follows:  For each pair
$X,Y$ of nonconsecutive edges of 
$K$, compute the minimum distance between the segments
$MD(X,Y)$, define $U_{\md}(X,Y)=\frac{\text{length}(X) 
\,\cdot\,  \text{length}(Y) }{{[MD(X,Y)]^{2}} }$, and sum:
$$U'_{md}(P) = \sum_{\text{all edges } X } \;\;\sum_{Y \neq X \text{
or adjacent}}U_{\md}(X,Y)\;\;.$$ 
This version of $\umd$ counts each
edge-pair twice, analogous to a double integral over (most of)
$K\times K$.  We write $U'_{\md}$ just to
distinguish from the original version \cite{Si2} that counted
each pair once.  To consider knots with varying numbers of
segments, we regularize by subtracting the energy
associated to a standard regular $n$-gon \cite{Si3} (or see \cite{Si5}).
Note that $U'_{\md}$ is scale invariant, so we can use any regular
$n$-gon and get the same number.  We define
$$E_{md}(P) = \upmd(P) - \upmd(\text{regular }n\text{-gon})\;\;.$$ 

The energy $\umd$ has been implemented  in several software
systems \cite{Wu,HuntKED,Sc}, and studied in 
\cite{MilRaw00,KauffHuang96,KauffHuang98}.

We shall show that for suitable polygonal approximations $P$
of a smooth curve $K$, 
$E_{md}(P) \approx E_{0}(K)$.
While the M\"obius Energy is defined for $C^{1,1}$ curves, our
proof requires that the knot $K$ be $C^2$ smooth.

\begin{figure}
\begin{center}
\includegraphics{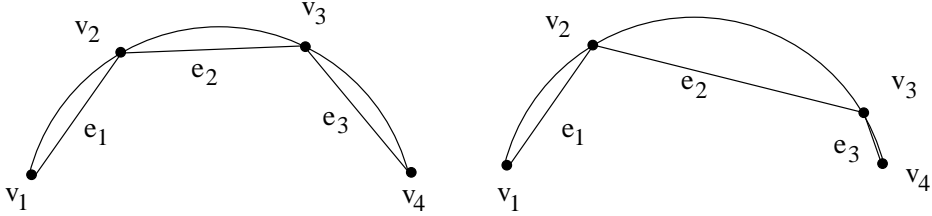}
\end{center}
\caption{Extreme differences in edge lengths causes $E_{md}$ to become large
despite having the polygon close to the smooth curve.  In such a case,
$E_{md}$ will be much larger than $E_0$.}
\vspace*{1em}
\label{unequal}
\end{figure}

In order for the approximation to work, we need to be careful
about what it means to say, ``the polygon P is a close
approximation of
$K$''.  First, we need to prevent extreme changes in the edge
lengths of
$P$ (see Figure \ref{unequal}). Suppose $P$ is a 
polygon closely inscribed in $K$.  We can slide vertex
$v_3$ along $K$ towards vertex $v_4$, making edge $e_3$
arbitrarily short, making edge $e_2$ longer, and keeping the
other edges of $P$ fixed.  This will make the contribution of
the edge-pair $(e_2,e_4)$ to
$E_{md}$ arbitrarily {\em large}.  Thus, we can make  polygons
$P'$ that also seem like close  approximations of $K$, yet
$E_{md}(P') >> E_{0}(K).$

To prevent this problem, we have to limit the variation in edge
lengths of polygons inscribed in $K$; we do this by having the vertices
equally spaced in arc-length along $K$. 
One can modify our arguments to handle other tractable approximating
polygons, e.g.~equal edge lengths or ``equal time'' subdivisions
of a regularly parameterized curve.

\begin{figure}
\begin{center}
\includegraphics{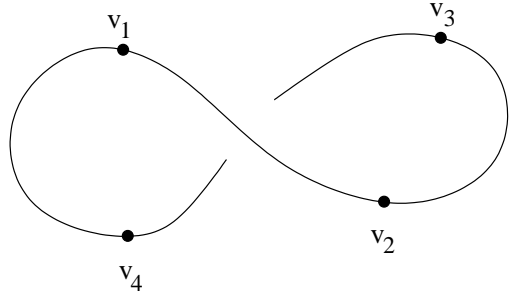}
\end{center}
\caption{Portions of the smooth curve can be arbitrarily close, causing
$E_0$ to be very large while $E_{md}$ of the inscribed polygon remains
fixed.  In such a case, $E_0$ is much larger than $E_{md}$.}
\label{polyhigh}
\end{figure}

Conversely, we can find situations where
$E_{md}(P) << E_{0}(K)$. Let $P$ be the polygon
(not drawn) $\langle v_1,v_2,v_3,v_4,v_1\rangle$ in Figure \ref{polyhigh}.
We construct the quadrilateral so that the arcs between
consecutive vertices are of equal length. Keeping the vertices
fixed, deform the arc $\widehat{v_1v_2}$ and the arc
$\widehat{v_3v_4}$ slightly so they get arbitrarily close
to intersecting (where one crosses over the other in the
figure).  This makes
$E_{0}(K)-E_{md}(P)$ arbitrarily large.  This problem is detected by the
fact that
$r(K)$ decreases to $0$, since normal disks of smaller and smaller
radii will intersect.  This is why the error bound in Theorem
\ref{explicit} must take into account the geometric quantity
$r(K)$.

\begin{figure}
\begin{center}
\includegraphics{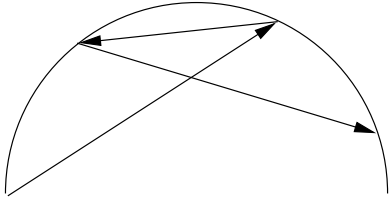}
\end{center}
\caption{We require that the order of the vertices of the inscribed
polygon coincides with an orientation of the smooth knot.}
\label{flipflop}
\end{figure}

To avoid situations as in Figure \ref{flipflop}, we
assume the phrase ``inscribed polygon'' means the vertices of
$P$ occur in the same order as they occur along $K$. 

Finally, note that the regularizations play an essential role, making 
the proof more delicate than may be at first evident.  See the discussion
in Section \ref{outlineproof}.

\section{Statement of main results and
notation}\label{MainResultsAndNotation}

We shall use the following  notation throughout the paper.
\begin{itemize}
\item $K$ is a $C^2$ smooth simple closed curve in
$\mathbb{R}^3$.
\item $\ell(K)$ is the total arc-length of $K$. 
\item $r(K)$ is the thickness radius of $K$. 
\item $E_L(K) = \frac{\ell(K)}{r(K)}$ is the Rope-Length.
\item $\delta=\frac{\ell(K)}{n}$ is the mesh size of the inscribed polygon.
\item $K$ is subdivided into $n$ arcs of equal length.
$\delta=\frac{\ell(K)}{n}$.
\item $v_1, \dots , v_n$ are the subdivision points along $K$. 
\item $P_n$ is the polygon formed by connecting the points
$v_1, \dots , v_n, v_1$ in order.
\item $\arc{x,y}$ is the length of the shorter of the two arcs
of $K$ connecting $x$ and $y$.
\item $|e|$ is the length of the line segment $e$.
\end{itemize}

Additional notation used in the proofs is listed at the
beginning of Section \ref{outlineproof}.

\begin{thm}[Error Bound] For any smooth
knot $K$, if $P_n$ are inscribed polygons as above
and $n$ is large enough that
$n>E_{L}(K)$,  i.e.~$\delta<r(K)$, then
$$|E_{0}(K)-E_{md}(P_n)|\leq \Phi(n,E_{L}(K))\;,$$
where $\Phi$ is a linear combination (see final page of 
the paper) of six fractions of the form $\frac{E_L(K)^a}{n^b}$ for
various $a>0$ and $b>0$.  By combining some terms, we can
take $\Phi = 550\frac{E_L(K)^{5/4}}{n^{1/4}}+ 10\frac{E_L(K)^4}{n}$.
For very large $n$ ($n>E_L(K)^{11/3}$), we can use
$\Phi = 560\frac{E_L(K)^{5/4}}{n^{1/4}}$.
\label{explicit}
\end{thm}

\noindent {\it Remark.}  There are other ways to write
this scale-invariant error bound, using
the identity
$$\frac{E_L(K)}{n} = \frac{\delta}{r(K)}\;.$$

From Theorem \ref{explicit}, we have immediately:

\begin{thm}[Approximation Theorem] For any smooth
knot $K$, if $P_n$ are inscribed polygons as above, then
as\; $n \to \infty$,
$E_{md}(P_n)\to E_{0}(K)$.
\label{approaches}
\end{thm}

\begin{pf}  The supporting
lemmas and the proof of the theorem occupy the rest of the paper. The
lemmas are in Section \ref{lemmas}. In Section \ref{outlineproof}, we
outline the proof and explain how the domains will be divided into
zones for which different analyses are needed. We give the analysis
for each zone in Section \ref{proofsforzones} and put them all
together in Section \ref{finalproof}.
\end{pf}

There are numerous coefficients in the calculations; we constantly
round up and pick the worst-case values, to keep the 
claims accurate and the numbers simple.

\section{The lemmas}
\label{lemmas}

In this section, we prove the lemmas needed for the proof of the main
theorem.

\subsection{Lemmas about the cosine function, also chords and arcs of circles}

\begin{lem} If $0< \phi\leq \pi$, then the following hold:
\begin{enumerate}
\item[(a)] $1-\frac{1}{2}\phi^2 \leq \cos \phi \leq 1-\frac{1}{2}\phi^2+
\frac{1}{24}\phi^4\,,$
\item[(b)]
$\phi^2-\frac{1}{12}\phi^4\leq 2-2\cos \phi \leq \phi^2\,,$
\item[(c)]
$1-\frac{1}{12}\phi^2\leq \frac{2-2\cos \phi}{\phi^2} 
\leq 1\,,$
\item[(d)]
$\frac{\phi^2}{2-2\cos \phi}
\leq 1+\frac{1}{2}\phi^2\,.$
\end{enumerate}
\label{allcosinelemmas}
\end{lem}
\begin{pf} For (a), consider the Taylor series for $\cos(\phi)$. 
Parts (b), (c), (d) follow immediately.
\end{pf}

\begin{lem} 
\mbox{}
\begin{enumerate}
\item[(a)]
On the unit circle  $C$, for any points $x,y$, 
$$\frac{1}{12} <    
\frac{1}{|x-y|^2} - \frac{1}{\mathrm{arc}(x,y)^2} \leq
\frac{1}{4}-\frac{1}{\pi^2}\,. $$
\item[(b)]
On a circle of radius $R$,
$$  \frac{1}{12}\;\frac{1}{R^2} <     
\frac{1}{|x-y|^2} - \frac{1}{\arc{x,y}^2} \leq 
\left(\frac{1}{4}-\frac{1}{\pi^2}\right)
\frac{1}{R^2}\,. $$
\end{enumerate}
\label{E0Integrand}
\end{lem}

\begin{pf}  Let $\phi$ be the angle ($\leq \pi$) between
points $x$ and
$y$ on the circle.  Since $C$ is the unit circle,
$\mathrm{arc}(x,y)=\phi$ and
$|x-y|^2=2-2\cos\phi$. The function
$$\frac{1}{2-2\cos\phi} -
\frac{1}{{\phi}^2}$$ is monotone, has a maximum at $\phi =
\pi$, and is bounded below by the limiting value as $\phi
\longrightarrow 0$.  Part (b) is similar.
\end{pf}

Next we want to compare the quantities  $\frac{1}{|x-y|^2}$ 
and  $\frac{1}{|X-Y|^2}$, where the points lie on circles of
different sizes.

\begin{lem}  Suppose $r<R$ are radii of circles and 
$0< a < \pi r$.  Construct any arcs of (the same) length
$ a $ on the two circles and let $x,y$ and $X,Y$ be the
endpoints of the two arcs.  Then  $$0 <
\frac{1}{|x-y|^2} - \frac{1}{|X-Y|^2} <  
\left(\frac{1}{4}-\frac{1}{\pi^2}\right)\,\frac{1}{r^2}\,.$$
\label{ChordsOnDifferentCircles}
\end{lem}

\begin{pf} Chord length is always
less than arc-length.  For a fixed arc-length, as the radius
gets larger, the chord length gets closer to the arc-length.  
Thus $|X-Y|>|x-y|$. 
On the other hand, applying Lemma
\ref{E0Integrand}(b) to each circle, we  have 
$$\frac{1}{|x-y|^2} - \frac{1}{|X-Y|^2} <  
\left(\frac{1}{4}-\frac{1}{\pi^2}\right)\;
\frac{1}{r^2} - \frac{1}{12}\;\frac{1}{R^2}.$$ 
 \end{pf}

\subsection{Lemmas about chords and arcs of general curves}

We rely a lot on Schur's Theorem.  Here is the version we need:

\begin{lem}  Let $K$ be a $C^2$ smooth curve in
$\mathbb{R}^3$ whose curvature everywhere is
$\leq$ some number $\kappa$.  Let $C$ be a circle of curvature
$\kappa$, i.e.~of radius $r =
\frac{1}{\kappa}$.  Let $x,y \in K, s,t \in C$ such that
$\mathrm{arc}(x,y)=\mathrm{arc}(s,t) \leq \pi r$.  Then the
chord distances satisfy
$$|x-y| \geq |s-t| \,.$$ 
When we write the chord length on $C$ in terms of the central
angle, this becomes
$$|x-y| \geq r \,\left( 2-2 \cos\left(\frac{\mathrm{arc}(s,t)}{r}   
\right)    
\right)^{1/2}
$$
\label{schursthm}
\end{lem}

\begin{pf} See Schur's Theorem in 
\cite{chern}.
\end{pf}

\begin{lem} Let  $K$ be a $C^2$ smooth curve in
$\mathbb{R}^3$, with minimum radius of curvature $r$. Suppose 
$x:[0,\pi r] \to \mathbb{R}^3$ is a unit speed
parameterization of an arc of $K$ of length $\pi r$.  Then the
function $|x(t)-x(0)|$ is monotone increasing.  That is to
say: As points move farther apart along the curve, they also move
farther apart in space, so long as the arc-distance is no greater
than $\pi r$.
\label{MonotoneArcs}
\end{lem}

\begin{pf} Let $f(t) = |x(t)-x(0)|^2 = (x(t)-x(0)) \cdot
(x(t)-x(0))$.  We claim $\frac{df}{dt}>0$ for $t\in (0,\pi
r)$.  The derivative
$\frac{df}{dt} = 2 (x(t)-x(0))\cdot x'(t)$.  Thus we 
need to show that this dot product is positive, for all points
$x(t)$ in the interior of the arc.  The proof uses the same
central idea as the proof of Schur's theorem.

We have
$$x(t)-x(0) = \int_0^t{x'(s)} \;ds \;,$$ 
so
$$(x(t)-x(0))  \cdot x'(t) = \int_0^t
x'(s)\cdot x'(t)  \;ds \;.$$ 
The dot product 
$x'(s)\cdot x'(t)$ is just the cosine of the angle $\leq \pi$
between the two velocity vectors.  This angle  is
measured by the length of the geodesic arc on the unit sphere
between the unit vectors
$x'(s)$ and  $x'(t)$. 
The trace of $x'(u)$, as $u$
runs from $s$ to $t$, is another path on the unit sphere
between the same vectors.  The length of that path gives an upper
bound for the length of the geodesic path.  Thus, since
$\left|x''(u)\right|\leq 1/r$ (recall $r =$ minimum radius
of curvature), 
$$\angle(x'(s),x'(t)) \leq \int_s^t{\left|x''(u)\right|} \;du
\;\leq \frac{(t-s)}{r}\;.$$

Since $0\leq s \leq t \leq \pi r$, and the cosine function is
decreasing on
$[0,\pi]$, we have 
$$\cos(\angle(x'(s),x'(t))) \geq \cos\frac{(t-s)}{r} \;.$$

Thus $$(x(t)-x(0))  \cdot  (x'(t)) \geq  \int_0^t
\cos\frac{(t-s)}{r}  \;ds \; = r \sin(t/r)\;.$$ 

For $0 < t < r \pi$, $\;\sin(t/r) > 0$.
\end{pf}

\begin{lem} Let  $K$ be a $C^2$ smooth curve in
$\mathbb{R}^3$, with minimum radius of curvature $r$.  Let
$x:[0,\ell(K)] \to \mathbb{R}^3$ be a unit speed
parameterization of $K$.   Suppose
$0\leq a<b<c<d\leq
\pi r$, so $x(a),x(b),x(c),x(d)$ are four points in order
along $K$, contained in an arc of total length $\leq \pi r$.  

Then the minimum spatial distance between line segments $\langle
x(a)x(b)\rangle$ and $\langle x(c)x(d)\rangle$ is realized at the
closest endpoints.  Taking into account Lemma \ref{MonotoneArcs},
this says,
$$ {MD}(\;\langle{x(a)x(b)}\rangle\;,\;\langle{x(c)x(d)}\rangle\;)
=
|x(c)-x(b)|\;.$$
\label{piarg}
\end{lem}
\begin{pf} Without loss of generality, rescale
the curve to have $r=1$.  Then the four points lie
in an arc of total length $\leq \pi$.

Let
$A$ denote the segment $\langle{x(a)x(b)}\rangle$ and $C$ the
segment 
$\langle{x(c)x(d)}\rangle$.  We shall show that for each
point $x\in A$, the point
$x(c)$ is the closest point of $C$ to $x$; so $x(c)$ is the
closest point of $C$ to $A$.  By a symmetric argument, 
the point $x(b)$
is the closest point of $A$ to $C$.

Fix a point $y \in C, y \neq x(c), x(d)$.  For any $x \in A$,
construct the directed line segment from $x$ to $y$.  We claim
that the vectors satisfy
$$(y-x) \cdot (x(d)-x(c)) > 0 \;.$$ If this dot product is
positive, then moving $y$ along $C$ closer to $x(d)$ will
increase the distance to $x$, and moving $y$ closer to $x(c)$
will decrease the distance to
$x$.  Thus $x(c)$ must be the closest point of $C$ to
$x$.

We now show that the above dot product is positive for each $x,y$.
It is convenient to think for a moment of fixing $y$ and varying $x$.
Let $P_y$ be the plane through $y$ perpendicular to $C$.  Rotate the
entire ensemble so that the vector $x(d)-x(c)$ points ``up''.  Then
the dot product inequality is equivalent to the assertion that the
entire line segment $A$ lies {\it below} $P_y$.  
It suffices to show that each vertex $x(a),
x(b)$ lies below $P_y$.

But in fact, if $x(a)$ and $x(b)$ lie
below $P_{x(c)}$, then they lie below $P_y$.  We have now
reduced the lemma to the following claim, an inequality that
involves only the given points on $K$.
The inequality is stated for parameter value
$a$,  and is identical for $b$.  If
$0
\leq a<c<d \leq \pi$, then 
$$ (x(c)-x(a))  \cdot  (x(d)-x(c)) > 0 \;.$$

The rest of the proof is similar to the proof of Lemma
\ref{MonotoneArcs} with some trigonometry at the end. We
first express the difference vectors as integrals of
derivatives,
$$(x(c)-x(a))  \cdot  (x(d)-x(c)) = \int_a^c\int_c^d
x'(s)\cdot\ x'(t) \;dt \;ds \;.$$ 

Since the cosine function is decreasing on $[0,\pi]$, we have 
$$\cos(\angle(x'(s),x'(t))) \geq \cos(t-s) \;.$$

As in the proof of Lemma \ref{MonotoneArcs},
$x'(s)\cdot x'(t) \geq \cos(t-s)$, so
$$(x(c)-x(a))  \cdot  (x(d)-x(c)) \geq  \int_a^c\int_c^d
\cos(t-s) \;dt \;ds \;.$$ 
The integral evaluates to 
$$\cos(d-c) - \cos(d-a) - 1 + \cos(c-a)\;,$$ which is positive.
\end{pf}

The previous two lemmas tell us that for arcs that are near each other
in arc-length along a curve, the minimum spatial distance between the
arcs is the same as the minimum distance between their inscribed
chords.  For more general pairs of arcs, the minimum distances usually
will not be equal, but they still are related.

\begin{lem} Suppose $\alpha,\beta$ are smooth arcs in 
$\mathbb{R}^3$, each  of
length
$\delta$, and each having radius of curvature everywhere $\geq
r\geq \delta$.  Let $e$ be the
chord joining the endpoints of $\alpha$ and $f$ the corresponding
chord for
$\beta$. 

(a) The maximum distance between
$\alpha$ and
$e$ (likewise between $\beta$ and $f$) is 
$\leq \frac{1}{\sqrt{48}}\frac{\delta^2}{r}.$

(b) If
 $\mathrm{MD}(\alpha,\beta)$ is the minimum spatial distance
between $\alpha$ and $\beta$, and $\mathrm{MD}(e,f)$ is the minimum
distance between the chords, then
$$|\mathrm {MD}(e,f)-\mathrm {MD}(\alpha,\beta)|\leq
\frac{\sqrt{3}\delta^2}{6r} \leq \frac{\sqrt{3}}{6}\,r\;.$$

\label{MaxDistArcToChord}
\end{lem}
\begin{pf} For part (a),
imagine the chord $e$ as a rod with ``string''
of length $\delta$ attached at either end, and ask, ``What configuration
allows the string to reach as far as possible from the
rod?''  The answer is  when the string is pulled out to form
two equal sides of an isosceles triangle, with the rod as the
base. The maximum distance that any point of $\alpha$ can be
from
$e$ is the altitude $h$ of this isosceles triangle, so
$h^2=\left(\frac{\delta}{2}\right)^2-\left(\frac{|e|}{2}\right)^2$.
Since $\delta \leq r$, in particular $\delta \leq \pi r$, we can
apply Schur's theorem: By Lemma
\ref{schursthm} and Lemma
\ref{allcosinelemmas}(b),
$$\frac{1}{4}|e|^2
\geq\frac{1}{4}\delta^2-\frac{1}{48}\frac{\delta^4}{r^2}
\;,$$ so
$$h^2\leq \frac{1}{48}\frac{\delta^4}{r^2}\,.$$

Part (b) follows from part (a), the triangle inequality, and the fact
that $\delta \leq r$.
\end{pf}

\subsection{Lemmas about the thickness of a curve}
\label{lemmasaboutthickness}

The first lemma is a characterization of the thickness
radius $r(K)$ in terms of curvature and the {\it critical
self-distance}.

Fix a  point
$x_0\in K$ and consider points
$y$ that start at $x_0$ and gradually move along $K$. A point
$y$ is  a critical point for the function $|y-x_0|$
when $y=x_0$ or when  $\langle xy \rangle \perp y'$.  We
define the {\it critical self-distance of $K$} (an idea
attributed  by J.~O'Hara to N.~Kuiper) to be
$$\mathrm{sd} (K) = \min\left\{\; |y-x|\;\;:\;\;x\not = y
\in K \;\mathrm{ and }\; \langle xy
\rangle \perp y' \;\right\}\;. $$

\begin{lem}  The thickness of a smooth knot is bounded by
the minimum radius of curvature and  half the  critical
self-distance.  In fact,
$$r(K)= \min \left\{ {\mathrm{MinRad}}(K), \frac12\; 
\mathrm{sd}(K)\right\}.$$
\label{lsdr}
\end{lem}
\begin{pf} See \cite{LSDR}.
\end{pf}

The next lemma is a consequence of Lemmas
\ref{schursthm} and \ref{lsdr}, and is proven in \cite{BRS}.
 
\begin{lem} Suppose $K$ is a smooth knot of thickness radius
$r(K)=r$.  For any $x,y \in K$ with
$\arc{x,y} \geq \pi r$, we must have $|y-x| \geq 2 r$.
\label{distancebound}
\end{lem}

\begin{lem}
Let $K$ be a $C^2$ smooth closed curve in $\mathbb{R}^3$, with minimum
radius of curvature $r$.  Let $C$ be a circle whose total arc-length
is the same as $K$, and $R$ be the radius of $C$.  Then $r \leq R$ and
(from Lemma \ref{lsdr}) the thickness radius $r(K) \leq R$.
\label{fenchellemma}\end{lem}
\begin{pf}
Since $r$ is the minimum radius of curvature of $K$, the
maximum curvature of $K$ is $\frac{1}{r}$, so the total
curvature of
$K$ is at most $\frac{\ell(K)}{r}$. On the other hand, by
Fenchel's theorem \cite{fenchel}, the total curvature of $K$
is at least $2 \pi$. Thus $2 \pi r \leq \ell(K) = 2 \pi R$.
\end{pf}

\begin{lem}
For any $C^2$ smooth closed curve $K$, $E_L(K)\geq 2\pi$.
\label{ropebiggerthan2pi}
\end{lem}
\begin{pf}
By Lemma \ref{lsdr}, the curvature of $K$ is everywhere $\leq 1/r(K)$.
Thus the total curvature of $K$ is $\leq \ell(K)/r(K)=E_L(K)$.
But the total curvature of a closed curve is $\geq 2\pi$.
\end{pf}

\section{Notation and outline of proof of Theorem \ref{explicit}}
\label{outlineproof}

\begin{figure}
\centering
\includegraphics{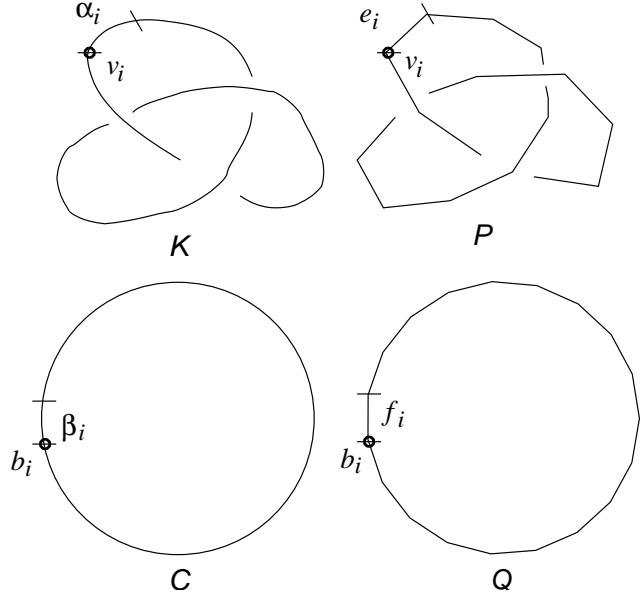}
\caption{The objects of study: smooth knot $K$ with arc $\alpha_i$ and
vertex $v_i$, inscribed polygon $P$ with vertex $v_i$, 
circle $C$ with arc $\beta_i$ and vertex $b_i$ corresponding to
$\alpha_i$ and $v_i$ respectively, and inscribed regular $n$-gon $Q$
with vertex $b_i$.}
\label{SetupNames}
\end{figure}

We have four objects of interest: the knot $K$, the circle $C$,
the inscribed $n$-gon $P$, and the regular  $n$-gon $Q$. In
the following list, refer to Figure \ref{SetupNames}

\begin{itemize}
\item $K$ is a $C^2$ smooth simple closed curve in
$\mathbb{R}^3$.
\item $C$ is a circle with total arc-length $\ell(C)=\ell(K)$.
\item $r(K)$ is the thickness radius of $K$. 
\item $K$ is subdivided into $n$ arcs of equal length
$\delta=\frac{\ell(K)}{n}$, and we are assuming 
$\delta<r(K)$ (so $n>E_L(K)$).
\item $v_1, \dots , v_n$ are the subdivision points along $K$.
\item $\alpha_i$ is the arc of $K$ with
endpoints $v_i$ and $v_{i+1}$.  We number the vertices modulo
$n$, so $\alpha_n$ is the arc from $v_n$ to $v_1$.
\item $R$ is the radius of $C$, so $R=\frac{\ell(K)}{2\pi}$.
\item $t \to x(t)$ is a unit speed parameterization of $K$
from $C$.
\item $b_1, \dots, b_n$ are evenly spaced points along $C$
such that $x(b_i)=v_i$.
\item $\beta_i$ is the arc of $C$ corresponding to $\alpha_i$.
\item $P$ is the polygon formed by connecting the points
$v_i$ in order.
\item $e_i$ is the edge of $P$ from $v_i$ to $v_{i+1}$, with length
denoted $|e_i|$.
\item $Q$ is the  regular polygon inscribed in $C$, with
vertices $b_1,\dots,b_n$.
\item $f_i$ is the edge of $Q$ with vertices $b_i, b_{i+1}$,
with length $|f_i|$.
\end{itemize}

Just to have all the important parameters specified in one
place, we also define two integers, $m$ and $p$, whose role
will be evident later in this section.
\begin{itemize}
\item $m = \lfloor \frac{\pi r(K)}{\delta} \rfloor$.  For a
vertex $v_i$, the vertices $v_i, v_{i+1}, \dots, v_{i+m}$
are a maximal list that lie in an arc of $K$ of length $\leq
\pi r(K)$.
\item $p = \lfloor m^{\frac{3}{4}}\rfloor$.  For a list of $m$
vertices as specified in the previous item, we will need to
distinguish an initial bunch from the rest.  It turns out
that the number we need to separate off should be some fractional power
of $m$ strictly greater than $1/2$, and we take
$3/4$ for simplicity.
\end{itemize}

We shall analyze the energies in terms of individual pairs of
arcs and/or edges. 
 
The energies are
\begin{align} E_{0}(K) &= \int_{x\in K}\int_{y\in K}
\frac{1}{|x-y|^2} - \frac{1}{|s-t|^2}\nonumber\\ &=
\sum_{i=1}^{n} \sum_{j=1}^{n} E_{0}(\alpha_i,\alpha_j)\;,
\label{esum}
\end{align} where $$E_{0}(\alpha_i,\alpha_j)=\int_{x\in
\alpha_i}\int_{y\in \alpha_j}
\frac{1}{|x-y|^2} - \frac{1}{|s-t|^2}\;,$$ 
and
\begin{align} E_{md}(P) &= U'_{md}(P)-U'_{md}(Q)\nonumber\\
&= \sum_{i=1}^{n} \;\sum_{j=1}^{n}
U_{md}(e_i,e_j)-U_{md}(f_i,f_j)\;\;\;(j \neq i-1,i,i+1)\;.
\label{usum}
\end{align} 
Sometimes we need to treat $E_0$ as the difference between two integrals, 
so we also define
\begin{equation*}
E(\alpha_i,\alpha_j) = \int_{x\in
\alpha_i}\int_{y\in \alpha_j}\frac{1}{|x-y|^2}\;,
\end{equation*}
and likewise for $E(\beta_i,\beta_j)$ for arcs on $C$.
As one
might expect, our overall plan is to show that the various
terms in the sum (\ref{esum}) are close to the  corresponding
terms in (\ref{usum}).  However, some terms in (\ref{esum})
have no corresponding term in (\ref{usum}); and even when
they do, there are different cases requiring different
analyses.  We shall, in fact, consider four kinds of pairs
$(i,j)$, bound each contribution to the error, and add them to
get a full error bound.

Here is a ``schematic diagram'' of our situation: We want to show that
something of the form $\int (W - X)$ is close to something of the form
$(Y-Z)$.  For the edge pairs where $E_0$ has a contribution and
$E_{md}$ is not defined, we show the $E_0$ contribution is small.  For
other edge pairs, we sometimes show that $\int (W - X)$ and $(Y-Z)$
each is small, and sometimes show that $|Y - \int W| $ and $|Z - \int
X|$ both are small. The analysis has to involve this kind of
complication because the unregularized polygon energy $U'_{md}(P)$ is
{\bf not} a good approximation of the divergent integral $\iint_{K
\times K} \frac{1}{|y-x|^2}\;$, that is $|Y - \int W| $ does not get
negligibly small for arc pairs (and their corresponding segment pairs)
that are extremely close together along $K$. Here is a simple example
to illustrate the difficulty: Consider two segments $A = [0,\epsilon]$
and $B = [2 \epsilon,3\epsilon] \subset \mathbb{R}$. Then $U_{md}(A,B)
= 1$.  On the other hand, $\int_{x\in A}\int_{y \in
B}\frac{1}{|y-x|^2}\;dy\;dx =\ln \frac{4}{3}\;.$ For segments close
together along the curves, we need to understand the regularizing
terms rather than show the two
energies are close to each other.

Following are the four types
of pairs (of indices $(i,j)$, edges or arcs) that determine
our four ``zones'' for separate analysis.  
The definitions are
symmetric, so $(i,j)$ and $(j,i)$ are of the same type.

\begin{enumerate}
\item {\it Adjacent Pairs}:  $j=i-1, i, i+1$\\ For these arc
pairs, we  bound $\sum_{i,j}{E_0(\alpha_i, \alpha_j)}$.
Since $\umd$ is only defined for non-adjacent edges,
there are no corresponding edge
pairs for these arc pairs.

\item {\it Near Pairs}: non-diagonal pairs $(i,j)$ for which
the arcs
$\alpha_i$ and $\alpha_j$ are contained in an arc of $K$ of
length
$\leq \pi r(K)$.\\ 
Within the Near Zone, we make an
additional distinction 
between ``Very Near'' and ``Moderately Near'': For each vertex $v_i$, let
$A$ be either of the arcs of $K$ starting at
$v_i$ and having length $\ell(A) = \pi r(K)$.   The vertices
contained in the arc $A$ are a sequence
$v_i, v_{i+1}, \hdots, v_{i+m}$ (for the other arc, we
count in the other direction).  The arcs contained in $A$ are
$\alpha_i,
\dots, \alpha_{i+m-1}$. The vertex
$v_{i+m}$ may or may not  be an endpoint of $A$. We
distinguish between the first
$m^{3/4}$ vertices and the rest. 
\begin{enumerate}
\item[A.] For $j=i+2, \hdots, i+p$,
we call $(i,j)$ a {\it very near} pair.\\ For such $(i,j)$, we
 bound $\sum_{i,j}\left(U_{md}(e_i,e_j)-U_{md}(f_i,f_j)\right)$
and 
$\sum_{i,j}{E_0(\alpha_i, \alpha_j)}$.
\item[B.] For $j = i+ p+ 1, \hdots, i+m-1$, we
call $(i,j)$ a {\it moderately near pair}.\\ For such
$(i,j)$, we shall bound
$\sum_{i,j}\left(E(\alpha_i,\alpha_j) -
U_{md}(e_i,e_j)\right)$ and
$\sum_{i,j}\left(E(\beta_i, \beta_j) -
U_{md}(f_i,f_j)\right).$
\end{enumerate}

\item \textit{Far Pairs}: The pairs $(i,j)$ that are neither
{\it adjacent} nor {\it near} are called 
{\it far}.\\   For such pairs, we shall also bound
$\sum_{i,j}{E(\alpha_i,\alpha_j) - U_{md}(e_i,e_j)}$ and
$\sum_{i,j}{E(\beta_i,\beta_j) - U_{md}(f_i,f_j)},$
but we need an argument different from the moderately near pairs.
\end{enumerate}

\begin{figure}
\begin{center}
\includegraphics[width=5.0in]{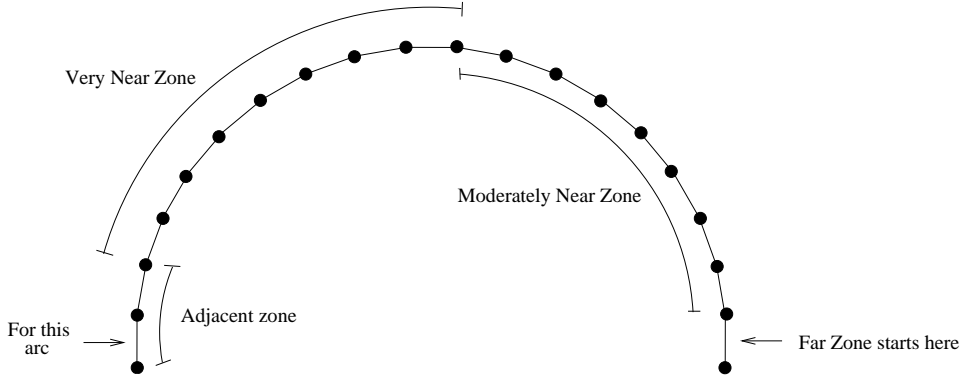}
\end{center}
\caption{The four types of zones on which we do our analysis.  
Note that this is just a schematic to show the arrangement of the
zones with respect to a fixed arc.}
\label{zones}
\end{figure}

\noindent See Figure \ref{zones} for an example of the zone
pairings where $m=17$.  We use the same terminology for corresponding
pairs of arcs in $C$; that is, if $(i,j)$ are far 
[resp.~adjacent, very near, moderately near] on $K$, then we
call them far [resp. adjacent, very near, moderately near]
on $C$.  

In the next section, we establish the explicit error bounds in each of
the different zones.  In Section \ref{finalproof}, we collect
all of the errors to determine the total error bound.

\section{Proofs for the different
zones}\label{proofsforzones}

\subsection{Bounds for $E_0$ in Adjacent and Very Near Zones}

We establish the error bound for
the combined contributions of the
Adjacent and Very Near Zones to
the M\"obius Energy.

\begin{prop} In the Adjacent and Very Near Zone,
$$\left|\sum_{i,j}{E_0(\alpha_i,
\alpha_j)}\right|<
1.06\,\frac{E_L(K)^{5/4}}{n^{1/4}}$$ 
\label{EoDiagAndVeryNear}
\end{prop}

\begin{pf} 
If $x,y$ are contained in diagonal or very near arcs, then
$\arc{x,y} \leq (p+1)\delta$.  Thus it suffices to bound
$$\left|\int_{x \in K}\int_{y \in K,\;
\mathrm{arc}(x,y)
\leq (p+1)\delta}\;\;
\frac{1}{|x-y| ^2}-\frac{1}{|s-t|
^2}\right|.$$

The calculation is independent  of
the choice of $x$, so we analyze
\begin{equation}
\left|2\,\ell(K)\,\int_{y=x}^{x+(p+1)\delta}
\frac{1}{|x-y|^2}-\frac{1}{|s-t|^2}\,dy\right|\,,
\label{eq:WantToBoundVeryNear}
\end{equation} 
where the limits of integration are meant to indicate that we
are integrating along an arc of $K$ of length $(p+1)\delta$
starting from $x$.

We are going to find upper and lower bounds for the integrand,
 observe that the upper bound is positive and the lower
bound is negative, and conclude that the magnitude of the
integrand is bounded by the difference between the upper and
lower bounds.  To simplify subsequent
expressions, let
$r$ denote $r(K)$ and $a$ denote $\arc{x,y}$.

Since
$s$ and
$t$ lie on a circle of radius $R$,
$$\frac{1}{|s-t|
^2}=\frac{1}{R^2(2-2\cos(a/R))}\,.$$

First we get the upper bound.  Since $\delta \leq r$, in
particular $m \geq 2$, we have $p <
m $ and $(p+1)\delta \leq \pi r$.  Thus we can apply Lemma
\ref{schursthm} to conclude
$$|x-y|^2 \geq r^2(2-2\cos(a/r))\;.$$

So we have
$$ \frac{1}{|x-y|^2}-\frac{1}{|s-t|^2}
\leq \frac{1}{r^2(2-2\cos(a/r))}
-\frac{1}{R^2(2-2\cos(a/R))}\,.
$$

By Lemma \ref{fenchellemma}, $r\leq R$. By Lemma
\ref{schursthm} applied to circles of different radii, or the
argument in Lemma \ref{ChordsOnDifferentCircles}, this upper
bound is nonnegative.

Now we get the lower bound. Since arc-length on any curve must
be at least as large as chord length, 
$$ |x-y|^2 \leq a^2\;.$$ 
 Thus
\begin{align}
\nonumber \frac{1}{|x-y|^2}-\frac{1}{|s-t|^2} &\geq
\frac{1}{a^2}-\frac{1}{|s-t|^2} \\ &= \nonumber
\frac{1}{a^2}-\frac{1}{R^2(2-2\cos(a/R))}\;,
\end{align} 
which is negative since chord length $<$ arc-length on a circle.

Taking the difference between the nonnegative upper bound
and the negative lower bound, we have
$$\left|\frac{1}{|x-y|^2}-
\frac{1}{|s-t|^2}\right| \leq
\frac{1}{r^2(2-2\cos(a/r))} - \frac{1}{a^2 \;.}
$$

So
$$(\ref{eq:WantToBoundVeryNear})
\leq 2\,\ell(K)\,\int_{0}^{(p+1)\delta}
\frac{1}{r^2(2-2\cos(a/r))}-\frac{1}{a^2}\,da\,,$$
where now we are just integrating a function of a real
variable.
Applying Lemma \ref{E0Integrand}(b), we  have
$$ (\ref{eq:WantToBoundVeryNear}) \leq 2\,\ell(K)\,
(p+1)\,\delta\,
\left(\frac{1}{4}-\frac{1}{\pi^2}\right)\,\frac{1}{r^2}\;.$$

Since $m \geq 3$, $p \geq 2$, so $(p+1) < 1.5 p$.  Combining
the constants, we have
\begin{equation*}
\nonumber (\ref{eq:WantToBoundVeryNear}) < 0.45\,
\frac{\ell(K) p
\delta}{r^2}
\leq  0.45\,\frac{\ell(K) (\frac{\pi
r}{\delta})^{3/4}
\delta}{r^2}
\leq 1.06\,
\frac{E_L(K)^{5/4}}{n^{1/4}}
\end{equation*}
 as desired.
\end{pf}

\subsection{Bound for $\emd$ in the Very Near Zone}

\begin{prop} In the Very Near Zone,
$$\left|E_{md}(P)\right| = |\umd'(P)-\umd'(Q)|<
2.76\, \frac{E_L(K)^{5/4}}{n^{1/4}}\,.$$
\label{polygonverynearzone}
\end{prop}

\begin{pf} 
$$\emd(\mathrm{very\ near})=
2\,\sum_{i=1}^{n}\sum_{j=i+2}^{i+p}
\frac{| e_i|\,| e_j|}{MD(e_i,e_j)^2}-
\frac{| f_i|\,| f_j|}{MD(f_i,f_j)^2}\,.$$ We shall bound the
inner sums uniformly in $i$, that is bound
\begin{equation}
\left| \sum_{k=1}^{p-1}\;
\frac{| e_i|\,| e_{i+k+1}|}{MD(e_i,e_{i+k+1})^2}-
\frac{| f_i|\,| f_{i+k+1}|}{MD(f_i,f_{i+k+1})^2}\right|
\label{eq:EmdVeryNearEq1} \end{equation}
 for arbitrary
$i$, and then multiply that bound by $2n$.  Here
$k=j-i-1$ is the number of edges separating the two edges.  As in
Proposition
\ref{EoDiagAndVeryNear}, we find a positive upper bound for each
difference term, and a negative lower bound; so the difference between
the upper and lower bounds is  a bound for the absolute value.

On the circle $C$ of radius $R$, the edge lengths are
$|f_i| = |f_j| =\sqrt{ R^2(2-2\cos(\delta/R))}\;$, and 
$MD(f_i,f_{i+k+1})=\sqrt{R^2(2-2\cos(k\delta/R))}$. So 
$$(\ref{eq:EmdVeryNearEq1}) \;=\;\left| \sum_{k=1}^{p-1}
\frac{| e_i|\,| e_{i+k+1}|}{MD(e_i,e_{i+k+1})^2}-
\frac{R^2(2-2\cos(\delta/R))}{R^2(2-2\cos(k\delta/R))}\right|\,.$$

To simplify subsequent expressions, let $r$ denote $r(K)$.
If we compare $K$ locally with a circle of radius $r$, 
Lemma \ref{piarg} and Lemma \ref{schursthm} say
$MD(e_i,e_{i+k+1})^2 \geq r^2(2-2\cos(k\delta/r))$.  The longest
an edge can be is the arc-length, so $(| e_i|\,| e_j|)
\leq \delta^2$.  Thus, an upper bound for each summand is
$$\mathrm{summand} \leq 
\frac{\delta^2}{r^2(2-2\cos(k\delta/r))}-
\frac{R^2(2-2\cos(\delta/R))}{R^2(2-2\cos(k\delta/R))}\,.$$
We claim this upper bound is positive. First, $\delta^2 >
R^2(2-2\cos(\delta/R))$ since arc-length (now on the big circle $C$)
 is always 
$>$ chord length.  Furthermore,
$r^2(2-2\cos(k\delta/r)) \leq R^2(2-2\cos(k\delta/R))$ by Lemma
\ref{ChordsOnDifferentCircles}.  

We next obtain a lower bound.  By Lemma \ref{lsdr}, $r \leq $ minimum
radius of curvature of $K$.  So we can apply Lemma \ref{schursthm} 
and Lemma \ref{piarg}
to  any points that lie in an arc of $K$ of length $\leq \pi r$.
By Lemma \ref{schursthm},  we have
$(| e_i|\,| e_{i+k+1}|) \geq r^2(2-2\cos(\delta/r))$.  
For the denominator, Lemma \ref{piarg} gives us that
$MD(e_i,e_{i+k+1}) = |v_{i+k+1}-v_{i+1}|$,
the distance between points of $K$ whose arc-distance is $k\delta$.
Since chord length $\leq$ arc-length, we thus have  
$MD(e_i,e_{i+k+1})^2\leq
(k\delta)^2$.  So a lower bound for the summand is
$$\mathrm{summand} \geq 
\frac{r^2(2-2\cos(\delta/r))}{k^2\delta^2}-
\frac{R^2(2-2\cos(\delta/R))}{R^2(2-2\cos(k\delta/R))\textsf{}}\,.$$ Comparing
numerators and denominators as we did for the upper bound, we see
that this lower bound is always negative.

Thus, we can bound the absolute value of the summand by the
difference between the upper and lower bounds:
\begin{align}
\left|\,\mathrm{summand}\,\right| &\leq
\frac{\delta^2}{r^2(2-2\cos(k\delta/r))}-
\frac{r^2(2-2\cos(\delta/r))}{k^2\delta^2}\nonumber\\
&=\frac{1}{k^2}\left(\frac{k^2\delta^2}{r^2(2-2\cos(k\delta/r))}-
\frac{r^2(2-2\cos(\delta/r))}{\delta^2}\right)\,.
\label{eq:UmdVeryNearEq2} 
\end{align}

We now appeal to our lemmas on cosines and chords.
To clarify how lemmas will be used, introduce angles 
$\theta=\delta/r$ and $\phi=k\delta/r$.  
 Thus, the bound (\ref{eq:UmdVeryNearEq2}) can be written
$$(\ref{eq:UmdVeryNearEq2}) =
\frac{1}{k^2}\left(\frac{\phi^2}{2-2\cos \phi}-
\frac{2-2\cos\theta}{\theta^2}\right)\,.$$

By Lemma \ref{allcosinelemmas}(d),     
 $\frac{\phi^2}{2-2\cos \phi}\leq
1+\frac{1}{2}\phi^2$.   By Lemma
\ref{allcosinelemmas}(c),  
$\frac{2-2\cos\theta}{\theta^2}\geq 1-\frac{1}{12}\theta^2$.
  Thus, 
$$|\mathrm{summand}|\leq
\frac{1}{k^2}\left(\frac{1}{2}\phi^2+\frac{1}{12}\theta^2\right).$$

We return to the original double sum and see that
\begin{align}
2n\,\left|\sum_{k=1}^{p-1}
\frac{| e_i|\,| e_{i+k+1}|}{MD(e_i,e_{i+k+1})^2}-
\frac{| f_i|\,| f_{i+k+1}|}{MD(f_i,f_{i+k+1})^2}\right| 
&\leq
2n\,\sum_{k=1}^{p-1}\frac{\frac{1}{2}\phi^2+\frac{1}{12}\theta^2}{k^2}
\nonumber\\ &=
2n\,\sum_{k=1}^{p-1}\frac{\frac{1}{2}k^2\theta^2+\frac{1}{12}\theta^2}{k^2}
\nonumber\\ &=
2n\,\frac{\delta^2}{r^2}\,\sum_{k=1}^{p-1}\left(\frac{1}{2}+
\frac{1}{12}\frac{1}{k^2}\right)\,. \nonumber
\\ &\leq 2 n\, \frac{\delta^2}{r^2} (p-1) \left(\frac{7}{12}\right)\nonumber
\\ &< \frac{7}{6}\, n\, \frac{\delta^2}{r^2} p \nonumber
\\ &< 2.76 \frac{E_L(K)^{5/4}}{n^{1/4}} 
\nonumber
\end{align}

\end{pf}

\noindent{\it Remark.} For the Very Near Zone, we could use $p\leq$
any fractional power $m^q$.  
It is in the Moderately Near Zone that we need $p > 1/2$.

\subsection{Bound for $|E_0(K)-\emd(P)|$ in the Moderately Near Zone}

In this section, we determine the error bounds in the Moderately
Near Zone for $|E(K)-\umd'(K)|$ and $|E(C)-\umd'(C)|$.  Recall that
the Moderately Near Zone consists of pairs $(i,j)$ where
$\alpha_i\,,\,\alpha_j$ 
[resp.~$\beta_i\,,\,\beta_j$] are contained in an arc of $K$ 
[resp.~$C$] of length $\pi\,r(K)$ but are separated by at least $p$ other
arcs; that is $k=j-i-1$ runs from $p$ to $(m-2)$.  The keys to the
analysis in this zone are: 
\begin{itemize}
\item The minimum distance between a given pair of arcs, or a given
pair of chords, is realized at the closest endpoints along the
curve.
\item That vertex-to-vertex distance is bounded away from zero by
Schur's theorem.
\end{itemize}

\begin{prop} \label{ModNearZone}
In the Moderately Near Zone,
$$|\teb| < 3.00\, \frac{E_L(K)^{11/4}}{n^{7/4}} + 542.84\,
\frac{E_L(K)^{3/2}}{n^{1/2}}\,.$$
\end{prop}
\begin{pf} 
As before, we use $r$ to abbreviate $r(K)$.  We first
analyze the error on $K$,
$$\left|2 \sum_{i=1}^{n}\sum_{k=p}^{m-2}\left( \frac{| e_i|\,| e_j
|} {MD(e_i,e_j)^2}-\int_{x\in \alpha_i}\int_{y\in \alpha_j}
\frac{1}{| x-y|^2}\,dy\,dx\right)\right|\,.$$
Note: The expressions seem more clear if we use both $k$ and $j$,
where $j=i+k+1$.

As in the previous case, the analysis is independent of $i$, so we
work with a general $i$ and multiply that bound by $n$.  To bound
the above sum of differences, we introduce a third term 
(larger than each of the two we are studying) and use
the triangle inequality.  

\noindent\textit{Claim 1}.\\   
$$2 n \sum_{k=p}^{m-2}\left| \frac{\delta^2} {MD(e_i,e_j)^2}-
\frac{| e_i|\,| e_j |} {MD(e_i,e_j)^2}\right| \leq
1.50\,\frac{E_L(K)^{11/4}}{n^{7/4}}\;.$$

\noindent \textit{Claim 2}.\\ 
$$2 n \sum_{k=p}^{m-2}\left|  \frac{\delta^2}
{MD(e_i,e_j)^2} 
- \int_{x\in \alpha_i}\int_{y\in
\alpha_j}
\frac{1}{| x-y|^2}\,dy\,dx \right| \leq
271.42\, \frac{E_L(K)^{3/2}}{n^{1/2}}\;.$$

\noindent {\it Proof of Claim 1.}
Since chord length $\leq$ arc-length, $|e_i||e_j| \leq \delta^2$.
So the summand without absolute value is non-negative, and any
upper bound will bound the absolute value.

Since we are still within the Near Zone, Lemma
\ref{schursthm} and Lemma \ref{allcosinelemmas}(b) give   
$$|e_i||e_j| \geq r^2(2-2 \cos(\delta/r)) \geq
\delta^2-\frac{1}{12}\frac{\delta^4}{r^2}\;.$$ 

Now consider the denominator.  By Lemmas \ref{piarg},
\ref{schursthm}, and \ref{allcosinelemmas}(b) 
\begin{align}
\nonumber MD(e_i,e_j)^2 = |v_j-v_{i+1}|^2  &\geq r^2(2-2
\cos(k\delta/r))
\\ \nonumber &\geq  k^2 \delta^2 - \frac{1}{12}\frac{k^4
\delta^4}{r^2}
\end{align}

Thus
$$2 n \sum_{k=p}^{m-2}\left( \frac{\delta^2} {MD(e_i,e_j)^2}-
\frac{| e_i|\,| e_j |} {MD(e_i,e_j)^2}\right) \leq
\frac{1}{6}\, n\,\delta^2 \,
\sum_{k=p}^{m-2}\frac{1}{k^2}\;\frac{1}{(r^2-\frac{1}{12}k^2
\delta^2)}\;.$$

We next bound this denominator away from $0$.  In the Near
Zone, $k \delta < \pi r$, so $r^2-\frac{1}{12}k^2
\delta^2 > r^2(1-\frac{1}{12}\pi^2)$, which gives
$$\frac{1}{6}\frac{1}{r^2(1-\frac{1}{12}\pi^2)} < 0.94\,
\frac{1}{r^2}\;.$$

We thus have
\begin{align}
2 n \sum_{k=p}^{m-2}\left( \frac{\delta^2} {MD(e_i,e_j)^2}-
\frac{| e_i|\,| e_j |} {MD(e_i,e_j)^2}\right) &\leq
0.94 \, n\,\frac{ \delta^2}{r^2} \,
\sum_{k=p}^{m-2}\frac{1}{k^2} \nonumber 
\\ &< 0.94 \, n\,\frac{
\delta^2}{r^2} \,
\sum_{k=p}^{\infty}\frac{1}{k^2} \nonumber \\
&< 0.94\, n\, \frac{\delta^2}{r^2}\,\frac{1}{p-1}\nonumber \;.
\end{align}

We want to bound $\frac{1}{p-1}$ in terms of $\delta$ and $r$.  Recall
that $p = \lfloor m^{3/4}\rfloor$ and $m = \lfloor 
\frac{\pi r}{\delta}\rfloor$.
Since $\delta < r$, we have $m \geq 3$ and $p \geq 2$.  So
$\frac{1}{p-1} \leq \frac{3}{p+1}$ and
$\frac{1}{m}\leq \frac{4}{3}\frac{1}{m+1}$.  Thus,
\begin{equation*}
\frac{1}{p-1} \leq 
\frac{3}{p+1} 
< \frac{3}{m^{3/4}}
 \leq 3\left(\frac{4}{3}\right)^{3/4}\frac{1}{(m+1)^{3/4}}
<  3.73\left(\frac{\delta}{\pi r}\right)^{3/4}
< 1.59\left(\frac{\delta}{r}\right)^{3/4}\,.
\end{equation*}

Then 
\begin{align}
2 n \sum_{k=p}^{m-2}\left( \frac{\delta^2} {MD(e_i,e_j)^2}-
\frac{| e_i|\,| e_j |} {MD(e_i,e_j)^2}\right) &< 
(0.94)(1.59) \, n\,\frac{
\delta^2}{r^2}\left(\frac{\delta}{r}\right)^{3/4} \nonumber
\\ &< 1.50\, \frac{E_L(K)^{11/4}}{n^{7/4}} \;. \nonumber
\end{align}

\noindent This completes the proof of Claim 1.

\bigskip\noindent {\it Proof of Claim 2.}

We need to bound 

\begin{equation}
\label{eq:Claim2ToBound}
2 n \sum_{k=p}^{m-2}\left|
\frac{\delta^2}{MD(e_i,e_j)^2}
- \int_{\alpha_i}\int_{\alpha_j}\;\frac{1}{|x-y|^2}\right|\;.
\end{equation}

By Lemma \ref{piarg}, $MD(e_i,e_j) = MD(\alpha_i,\alpha_j) =
|v_j-v_{i+1}|$.  Since the arcs have length $\delta$, we know that the
summands without absolute value are nonnegative; so, as in Claim 1, we
bound the absolute value by finding an upper bound.  We are dealing
with something that looks like a Riemann Sum upper estimate of a
finite integral.  But as $n$ increases, we are changing the domain,
not just subdividing the same set
and we want to control the size of the error, not just say it goes
to zero as $n \rightarrow \infty$.  This is where we use the
choice of $p$ as a fractional power $m^q$ where $q$ is strictly larger
than $1/2$.

For brevity, let
$md$ denote 
$MD(\alpha_i,\alpha_j)=|v_j-v_{i+1}|$, where $\arc{v_{i+1},v_j}=k\delta$.  
Since
$$|x-y| \leq md+2\delta\;,$$
we have
\begin{equation*}
\frac{\delta^2}{md^2} -
\int_{\alpha_i}\int_{\alpha_j}\;\frac{1}{|x-y|^2}
\nonumber 
\leq
\frac{\delta^2}{md^2} -\frac{\delta^2}{(md+2\delta)^2}
< 4 \delta^3 \frac{1}{md^3 }\,.
\end{equation*}

\noindent Thus
$$ (\ref{eq:Claim2ToBound})
\;\;\leq\;\; 8 n \delta^3 \sum_{k=p}^{m-2} \frac{1}{md^3}\;.
$$

As before, by Lemmas \ref{piarg},
\ref{schursthm}, and \ref{allcosinelemmas}(b),
since $\arc{v_{i+1},v_j}=k\delta$,
\begin{align}
\nonumber md^2 = |v_j-v_{i+1}|^2  &\geq r^2(2-2
\cos(k\delta/r))
\\ \nonumber &\geq  k^2 \delta^2 - \frac{1}{12}\frac{k^4
\delta^4}{r^2}
\\ \nonumber &= k^2 \delta^2
\left(1-\frac{1}{12}\frac{k^2\delta^2}{r^2}\right)
\\ \nonumber &\geq k^2 \delta^2
\left(1-\frac{1}{12}\pi^2\right) \text{since $k\delta \leq \pi r$.}
\end{align}

Thus, $\md\geq 0.42\,k\delta$, so
$$\frac{1}{md^3} < 13.42\, \frac{1}{k^3\delta^3} \;.$$ 

 With the above observation, the $\delta^3$'s
cancel and we have
\begin{align}
(\ref{eq:Claim2ToBound})
&< 107.36\, n  \sum_{k=p}^{m-2} \frac{1}{k^3} \nonumber \\
&< 107.36\, n  \sum_{k=p}^{\infty} \frac{1}{k^3} \nonumber \\
&< 107.36\, n\frac{1}{(p-1)^2}\; . \nonumber
\end{align}

We showed in the proof of the prior claim that
$\frac{1}{p-1}< 1.59 \left(\frac{\delta}{r}\right)^{3/4}$.  Thus,
\begin{equation*}
(\ref{eq:Claim2ToBound})
\leq 107.36\,n\frac{1}{(p-1)^2}
\leq 271.42\, \frac{E_L(K)^{3/2}}{n^{1/2}}
\end{equation*}

Note that in the above analysis
the exponent $3/4$ needs to be strictly greater than $1/2$, so that
when we double it, the power of $n$ in the denominator will more
than cancel the leading factor $n$.

We now need to bound the contribution from $C$, that is
$|E(C)-\umd'(C)|$.  The radius of $C$, $R$, is the thickness
radius $r(C)$.  Also, we know from Lemma \ref{ChordsOnDifferentCircles}
that $R\geq r$.  So if arcs $\alpha_i$, $\alpha_j$ of $K$ are near,
then the corresponding arcs $\beta_i$, $\beta_j$ lie within an arc
of $C$ of length $\leq \pi R$.  Thus, the various steps in our analysis
of $K$ can be carried out on $C$.  We could obtain sharper bounds for $C$,
but we will settle for the same bound since they dominate anyway.

For Claim 1, we have
$$|f_i|\,|f_j| = R^2(2-2\cos(\delta/R)) \geq \delta^2-
\frac{1}{12}{\delta^4}{R^4} \geq \delta^2-\frac{1}{12}\frac{\delta^4}{r^2}\,,$$
and 
$$MD(f_i,f_j)^2\geq k^2\delta^2 - \frac{1}{12}\frac{k^4\delta^4}{R^2}
\geq k^2\delta^2-\frac{1}{12}\frac{k^4\delta^4}{r^2}\,,$$
exactly as for $K$.  Now continue the proof of Claim 1 verbatim.

For Claim 2,
\begin{align*}
md(f_i,f_j)^2
&=R^2(2-2\cos(k\delta/r))\\
&\geq k^2\delta^2\left(1-\frac{1}{12}\frac{k^2\delta^2}{R^2}\right)\\
&\geq k^2\delta^2\left(1-\frac{1}{12}\frac{k^2\delta^2}{r^2}\right)\,,
\end{align*}
and the rest follows verbatim.

Thus, our final bound for the total error in this zone is just
double the values obtained in Claims 1 and 2.
\end{pf}

\subsection{Bounds for $|E_0(K)-\emd(P)|$ in the Far Zone}

As before, we use $r$ to abbreviate $r(K)$.  In the Near Zones, we
just needed a value for $r\leq$ minimum radius of curvature
of $K$.  But in the Far Zone, we need both aspects of the thickness radius.

The argument here is somewhat similar to the Moderately Near Zone,
but we control the denominators in a different way.  In each
situation, we need to know that spatial distances between points are
bounded away from zero in some way depending on their arc-length
distances along $K$.  
For the Far Zone, we use the fact that thickness
controls critical self-distance, in particular   Lemma
\ref{distancebound}, together with local analysis (Lemma
\ref{MaxDistArcToChord}), to relate
chord-chord distances to arc-arc distances. 
Also, we continue to use the hypothesis $\delta \leq r$.

\noindent\textit{Remark on notation heuristics}. In the following
paragraphs and Lemma \ref{FarChordsMin}, think of 
$(\alpha,\beta)$ as $(\alpha_i,\alpha_j)$ and $(e,f)$ as
$(e_i,e_j)$.

\begin{lem} 
Suppose $(\alpha,\beta)$ is a pair of far arcs (on
$K$ or on $C$), with
$(e,f)$ the inscribed chords joining their endpoints.  Then
$$\mathrm{md}(\alpha,\beta) > 1.08\,r\;,$$
and
$$\md(e,f)> 0.79 \, r\;.$$
\label{FarChordsMin} 
\end{lem}

\begin{pf}
We analyze $K$, and note that the same bound will work for $C$ since
$r\leq R$.
We establish the lower bound for arcs, then use that to bound
the distance for chords.  If the minimum distance between a pair of arcs
is realized at points that are interior to one or both arcs,
then we are dealing with singly- or doubly-critical pairs of
points, so, by Lemma \ref{lsdr}, $\md(\alpha,\beta) \geq 2r$.  Thus
we just need to bound the end-point distances.  
Let $\alpha_0$ and $\alpha_1$ be the endpoints of the arc $\alpha$ and
$\beta_0$ and $\beta_1$ the endpoints of the arc $\beta$.  Choose the
labels so
that $\alpha_1$ and $\beta_0$ are the points which are closest with
respect to arc-length.  In the worst case, the arc-length from
$\alpha_0$ to $\beta_1$ is $\geq \pi r$, but the arc-lengths of the arcs
$\widehat{\alpha_0\beta_0}$, $\widehat{\alpha_1\beta_0}$, and 
$\widehat{\alpha_1\beta_1}$ are less than $\pi r$.  
In such a case, we have the following situation:

\begin{itemize}
\item $|\alpha_0-\beta_0|$ \\
$\pi\,r \geq {\rm
arc}(\alpha_0,\beta_0)
\geq
\pi r-\delta
\implies |\alpha_0 -
\beta_0|^2 \geq r^2(2-2\cos(\pi-1))$ by Lemma \ref{schursthm}, and
the fact that $\delta \leq r$.  So $|\alpha_0-\beta_0|> 1.75 r$.

\item $|\alpha_0-\beta_1|$ \\
${\rm arc}(\alpha_0,\beta_1) \geq \pi r \implies
|\alpha_0 -
\beta_1| \geq 2r$ by Lemma \ref{distancebound}.

\item $|\alpha_1-\beta_1|$ \\
same bound as $|\alpha_0-\beta_0|$ .

\item $|\alpha_1-\beta_0|$ \\
${\rm arc}(\alpha_1,\beta_0) \geq \pi r-2\delta \implies
{\rm arc}(\alpha_1,\beta_0) > (\pi-2)r$, since $\delta \leq r$.  Thus, by
Lemma \ref{schursthm},
\\ $|\alpha_1-\beta_0|^2 \geq r^2(2-2 \cos(\pi-2)) 
\implies |\alpha_1-\beta_0| > 1.08r$.

\end{itemize}

In other scenarios, the arc pair $(\alpha,\beta)$ yields three of
the above four cases, but we lose the smallest.  For ``most'' arc pairs
$(\alpha,\beta)$, we have all point-to-point distances at least
$2r$. 

We now obtain the lower bound on chord-to-chord distances using 
Lemma \ref{MaxDistArcToChord}:
$$\md(e,f) \geq   \md(\alpha,\beta) - \frac{\sqrt{3}}{6}\,r > 
1.08\,r -\left(\frac{\sqrt{3}}{6}\right) r > 0.79\, r\;.$$

\end{pf}

\begin{prop}
\label{FarZoneBound}
 The total error in the Far Zone is
bounded by
$$0.56\,\frac{E_L(K)^4}{n^2}+1.60\,\frac{E_L(K)^5}{n^2}+7.76\,\frac{E_L(K)^4}{n}\,.$$
\end{prop}
\begin{pf} We first analyze the error on $K$,
\begin{equation*}
2\sum_{i=1}^n\sum_{j=i+m}^{n}\left|
\frac{|e_i|\,|e_j|}{\md(e_i,e_j)^2}-
\int_{x\in \alpha_i}\int_{y\in \alpha_j}
\frac{1}{| x-y|^2}\,dy\,dx\,\right|\,.
\end{equation*}

We do this in three steps: Compare $\frac{|e_i|\,|e_j|}{md(e_i,e_j)^2}$
to $\frac{\delta^2}{md(e_i,e_j)^2}$, that to
$\frac{\delta^2}{md(\alpha_i,\alpha_j)^2}$, and that to
$\iint\frac{1}{|x-y|^2}$.
After we do each step for $K$, we double that to include
the contribution from $C$.  Note 
$\frac{\delta^2}{md(e_i,e_j)^2} = 
\int_{x\in\alpha_i}\int_{y\in\alpha_j}\frac{1}{md(e_i,e_j)^2}\,dy\,dx$
and similarly for $\frac{\delta^2}{md(\alpha_i,\alpha_j)^2}$.

\noindent\textit{Claim 1:}\\
\begin{equation*}
2\sum_{i=1}^{n}\sum_{j=i+m}^n \left|
\frac{\delta^2}{md(e_i,e_j)^2}
-\frac{|e_i|\,|e_j|}{md(e_i,e_j)^2}\right|
\leq 0.28\,\frac{E_L(K)^4}{n^2}\,.
\end{equation*}

\noindent\textit{Claim 2}.
\begin{equation}
2\sum_{i=1}^{n}\sum_{j=i+m}^n \int_{x\in\alpha_i}\int_{y\in\alpha_j}
\left| \frac{1}{md(e_i,e_j)^2} -
\frac{1}{md(\alpha_i,\alpha_j)^2}\right|\,dy\,dx\,\leq
0.80\,\frac{E_L(K)^5}{m^2}.
\label{eq:mdedgearc}
\end{equation}

\noindent\textit{Claim 3}.
\begin{equation}
2\sum_{i=1}^{n}\sum_{j=i+m}^n \int_{x\in\alpha_i}\int_{y\in\alpha_j}
\left|\frac{1}{md(\alpha_i,\alpha_j)^2} -
\frac{1}{| x-y|^2}\right|\,dy\,dx\,\leq 3.88\,\frac{E_L(K)^4}{n}\;.
\label{lastone}
\end{equation}

\noindent\textit{Proof of Claim 1}.

Since arc-length $\geq$ chord length, each summand is nonnegative
without taking the absolute value, so we just need to bound
the terms from above.  By Lemma
\ref{schursthm} and Lemma
\ref{allcosinelemmas}(b), 
$\delta^2-\frac{1}{12}\frac{\delta^4}{r^2}\leq
|e_i|,|e_j|$. Thus,
$$\frac{\delta^2}{\md(e_i,e_j)^2}
-\frac{|e_i|\,|e_j|}{\md(e_i,e_j)^2}
\leq \frac{1}{12}\;\frac{\delta^4}{r^2 \;\md(e_i,e_j)^2}\,.$$ 
But Lemma \ref{FarChordsMin} gives us that
$\md(e_i,e_j)^2 >  (0.79)^2\,r^2$, so
$$\frac{\delta^2}{\md(e_i,e_j)^2}
-\frac{|e_i|\,|e_j|}{\md(e_i,e_j)^2}
< 0.14\,\frac{\delta^4}{r^4}\,.$$ 
Multiplying by $2n^2$ gives
\begin{equation*} 
2\sum_{i=1}^{n}\sum_{j=m}^n \left(
\frac{\delta^2}{\md(e_i,e_j)^2}
-\frac{|e_i|\,|e_j|}{\md(e_i,e_j)^2}\right) 
< 0.28\,\frac{n^2\delta^4}{r^4}
= 0.28\,\frac{E_L(K)^4}{n^2}
\end{equation*}

\noindent\textit{Proof of Claim 2}. 

The sum \eqref{eq:mdedgearc} is bounded by
$(2n^2\delta^2)\mathrm{(worst\ error\ in\ integrands)}.$
We will use Lemma \ref{MaxDistArcToChord}(b) to bound that.   
To make the algebra more evident, let $\epsilon = md(e_i,e_j)$
and $\gamma = md(\alpha_i,\alpha_j)$.
The term we wish to bound is
$$\left|\frac{1}{\epsilon^2}-\frac{1}{\gamma^2}\right|=
\left|\frac{\gamma^2-\epsilon^2}{\epsilon^2\gamma^2}\right|
< 1.38\,\frac{|\gamma-\epsilon|\,(\gamma+\epsilon)}{r^4}
\leq 0.40\,\frac{(\gamma+\epsilon)\delta^2}{r^5}\,,$$
since $\epsilon > 0.79 r$ and $\gamma>1.08r$ by Lemma \ref{FarChordsMin},
and $|\gamma-\epsilon|\leq \frac{\sqrt{3}}{6}\frac{\delta^2}{r}$ by
Lemma \ref{MaxDistArcToChord}(b).  

Now $\epsilon$, $\gamma$ are minimum distances between sets that include
points of $K$, so
$\epsilon,\gamma\leq \ell(K)/2$ and 
$\gamma+\epsilon\leq \ell(K)$. Thus,
$$\left|\frac{1}{\epsilon^2}-\frac{1}{\gamma^2}\right| \leq
\frac{0.40\,\ell(K)\delta^2}{r^5}\,.$$ 
Multiplying by $2n^2\delta^2$, we get
\begin{align} 2\sum_{i=1}^{n}\sum_{j=i+m}^n 
\int_{x\in\alpha_i}\int_{y\in\alpha_j}
\left|\frac{1}{md(e_i,e_j)^2} -
\frac{1}{md(\alpha_i,\alpha_j)^2}\right|\,dy\,dx &\leq 0.80\,
\frac{n^2\delta^4\ell(K)}{r^5}\nonumber\\ &= 0.80\,
\frac{E_L(K)^5}{n^2}\nonumber\,.
\end{align}

\noindent\textit{Proof of Claim 3}.

The sum \eqref{lastone} is bounded by
$(2n^2\delta^2)\mathrm{(worst\ error\ in\ integrand)}.$

Let $\gamma$ denote $md(\alpha_i,\alpha_j)$.
So for particular $x$, $y$ on $\alpha_i$ and $\alpha_j$, we have 
$|x-y|= \gamma+t$ for some
$0\leq t\leq 2\delta$.
The largest error is then
$$\frac{1}{\gamma^2}-\frac{1}{(\gamma+t)^2}=
\frac{t(2\gamma+t)}{\gamma^2(\gamma+t)^2}
\leq \frac{t(2\gamma+t)}{\gamma^4}
<\frac{t(2\gamma+t)}{(1.08)^4r^4}
\,,$$ 
since $\gamma\geq 1.08r$ by Lemma \ref{FarChordsMin}.

Now $t\leq 2\delta$ and $\gamma\leq \ell(K)/2$.  Thus,
\begin{align*}
\frac{t(2\gamma+t)}{1.08^4r^4}
&\leq \frac{2\delta(\ell(K)+2\delta)}{1.08^4r^4}\\
&= \frac{2}{1.08^4}\frac{\delta(n\delta+2\delta)}{r^4}\\
&= \frac{2}{1.08^4}\frac{\delta^2(n+2)}{r^4}\\
&\leq \frac{2}{1.08^4}\frac{\delta^2}{r^4}\frac{(2\pi+2) n}{2\pi}
\text{ since }n>\elk\geq 2\pi \\
&< 1.94\,\frac{\delta^2n}{r^4}\,.
\end{align*}

Thus, 
\begin{align*}
2\sum_{i=1}^{n}\sum_{j=i+m}^n \int_{x\in \alpha_i}
\int_{y\in \alpha_j} \left|\frac{1}{md(\alpha_i,\alpha_j)^2} -
\frac{1}{| x-y|^2}\right|\,dy\,dx
&\leq 3.88\, n^2\delta^2\frac{\delta^2n}{r^4}\\
&= 3.88\,\frac{E_L(K)^4}{n}
\end{align*}

\end{pf}

\section{Putting it all together}\label{finalproof}

Here we combine the bounds from the various zones.

From Propositions \ref{EoDiagAndVeryNear}, \ref{polygonverynearzone},
\ref{ModNearZone}, and \ref{FarZoneBound}, we have
\begin{align*}
|E_0(K)-\emd(P)| \leq&\; 3.82\,\elk\,\left(\frac{\elk}{n}\right)^{1/4}\\
&+ 3.00\,\elk\, \left(\frac{\elk}{n}\right)^{7/4}\\
&+ 542.84\,\elk\,\left(\frac{\elk}{n}\right)^{1/2}\\
&+ 0.56\,\elk^2\,\left(\frac{\elk}{n}\right)^2\\
&+ 1.60\,\elk^3\,\left(\frac{\elk}{n}\right)^2\\
&+ 7.76\,\elk^3\left(\frac{\elk}{n}\right)
\end{align*}

Since $E_L(K)\geq 2\pi>1$, and $n>\elk$, we see that certain terms dominate
others.  So,
$$|E_0(K)-\emd(P)| < 550\, \frac{\elk^{5/4}}{n^{1/4}} + 
10\, \frac{\elk^4}{n}\;.$$
If $n>\elk^{11/3}$, then the total error is less than 
$560\, \frac{\elk^{5/4}}{n^{1/4}}$.

This completes the proof of Theorem \ref{explicit}.

\section{Acknowledgments}

We thank J.~Sullivan for asserting and experimentally confirming
the correct regularization and Y.-Q.~Wu for
modifying {\it MING} to allow additional numerical confirmation.  We
also thank G.~Buck for helpful comments.

\bibliographystyle{elsart-num.bst}
\bibliography{JonsEricBib04-23-03}

\end{document}